\documentclass{amsart}
\usepackage{graphicx,,amssymb}
\usepackage{amsfonts}
\usepackage{amscd}
\usepackage{mathrsfs}
\usepackage[all]{xy}
\usepackage{overpic}

\newtheorem{thm}{Theorem}[section]
\newtheorem{lemma}[thm]{Lemma}
\newtheorem{prop}[thm]{Proposition}
\newtheorem{cor}[thm]{Corollary}

\newtheorem*{thm3}{Theorem 3.2}

\theoremstyle{definition}
\newtheorem{defn}[thm]{Definition}

\theoremstyle{remark}

\numberwithin{equation}{section}

\newcommand{\ad}{\text{ad}}

\newcommand{\tr}{\text{tr}}
\newcommand{\utA}{\underset{\sim}{A}}
\newcommand{\otA}{\overset{\sim}{A}}

\newcommand{\Ham}{\mathbb{H}}        

\newcommand{\CP}{\mbox{$\mathbb{CP}$}}        
\newcommand{\HP}{\mbox{$\mathbb{HP}$}}         
\newcommand{\RP}{\mbox{$\mathbb{RP}$}}         
\newcommand{\OP}{\mbox{$\mathbb{OP}$}}         

\DeclareMathOperator{\Hom}{Hom}

\begin{document}
\title{Rigidity of Rank-One Factors of Compact Symmetric Spaces}
\author{Andrew Clarke}
 \address{Centre des Math{\'e}matiques Laurent Schwartz, Ecole Polytechnique,  
            91128 Palaiseau Cedex, France}
    \email{clarke@math.polytechnique.fr}

\date{\today}

\maketitle

\begin{abstract}
We consider the decomposition of a compact-type symmetric space into a product of factors and show that the rank-one factors, when considered as totally geodesic submanifolds of the space, are isolated from inequivalent minimal submanifolds.
\end{abstract}


\section{Introduction}



Questions of isolation phenomena for minimal submanifolds have been posed for many years. Perhaps the most studied case is for minimal submanifolds of the sphere. 
 Lawson \cite{Lawson1}, Chern, do Carmo and Kobayashi \cite{CCK}, Barbosa \cite{Barbosa}, Fischer-Colbrie \cite{FC} and others studied minimal submanifolds of the sphere using a range of techniques and obtained existence and uniqueness results. An important part of this study was initiated by Simons \cite{Simons}, who used a rigidity-isolation result for minimal hypersurfaces of $S^n$ to show that a minimal cone in euclidean space constructed as a blow-up limit from the given minimal graph was over a totally geodesic subset in the sphere. This was an important part of his extension of the Bernstein theorem to dimensions up to $n=7$.

Minimal submanifolds of other specific geometric spaces have also been studied. Much work has been done on the classification of totally geodesic submanifolds of riemannian symmetric spaces. This is particularly tractable because of the identification of these spaces with Lie triple systems. This allows the representation theory and algebra of the ambient isometry group to be considered. Gluck, Morgan and Ziller \cite{GMZ} and Thi \cite{Thi} used the terminology of calibrations to study the homologically volume minimizing cycles in Grassmannians and Lie groups respectively. W.T. Hsiang and W.Y. Hsiang \cite{Hsiang} constructed non-totally geodesic minimal hypersurfaces diffeomorphic to spheres in a range of compact-type symmetric spaces. 

In the case at hand we consider the decomposition of a compact-type Riemannian symmetric space into irreducible components and consider the closed minimal submanifolds that are close, in a concrete sense, to the factors. If the factor has rank equal to one and has non-exceptional isometry group, any nearby closed minimal submanifold must be another factor in the decomposition. The nearness that we consider implies that the submanifold is the graph of a map from the factor to the other components so the result can also be thought of as an isolation-type statement for maps between symmetric spaces.

\begin{thm}  
Let $X$ be a symmetric space of compact type. Suppose that $X$ decomposes as $X=X_1\times X_2$ where $X_1$ is of rank-one, and is not the Cayley Plane. Then there is a $C^3$ neighbourhood of the standard embedding of $X_1$ as a factor such that any minimal immersion in this neighbourhood is as another factor in the decomposition.
\end{thm}

This result follows from Theorems \ref{thm-tgisol}, \ref{CPnRigid} and \ref{HPnrigid}.

This work was completed while the author was a student at SUNY Stony Brook. He would like to acknowledge the enormous generosity and guidance given by his advisor, Blaine Lawson.


\section{Preliminaries}

We recall the important calculation of Simons that showed that the second fundamental
 form of a minimal submanifold satisfies a second order elliptic equation. That is, we define the {\em Second Fundamental Form} $A$  of a submanifold $M\subseteq X$ to be 
\begin{eqnarray*}
A^{\nu}Y =-(\nabla_Y\nu)^T
\end{eqnarray*}
where $\nu$ is a normal vector to $M$ and $Y$ is tangent to $M$. A priori, $\nu$ must be defined locally but it is clear that this definition is independent of the extension. $A$ is in thus a section of the Riemannian vector bundle $\text{Hom}(N_M, S(M))$ of endomorphisms from the normal bundle to symmetric transformations of the tangent bundle. We say that $M\subseteq X$ is {\em minimal} if $A^\nu$ is trace-free for all normal vectors $\nu$. Simons \cite{Simons} was able to calculate the rough Laplacian  of $A$. He showed that if $M$ is minimal, $\nabla^2 A$ can be algebraically expressed only in terms of $A$ itself, together with the ambient curvature and its covariant derivative. That is,
\begin{eqnarray*}
\nabla^2 A= -A\circ \otA - \utA\circ A +R(A) +R^\prime.
\end{eqnarray*}
$\utA$ and $\otA$ are quadratic expressions in $A$ and $R^\prime\equiv 0$ if $X$ is locally symmetric. He also showed 
$\utA$ and $\otA$ universally and uniformly satisfy the inequality
\begin{eqnarray*}
\langle A\circ \otA+\utA\circ A,A\rangle\leq q\|A\|^4
\end{eqnarray*}
where $q=2-\frac{1}{\text{codim} M}$. With this inequality, and the assumption that $X$ is symmetric, we can largely overlook these terms from here. The exact definition of $\otA$ and $\utA$ are given in \cite{Simons} but we will only require this inequality.

The other term $R(A)$ is also a section of $\text{Hom}(N_M,S(M))$ and is given by 

\begin{eqnarray*}
\langle R(A)^WX,Y\rangle &=& 
\sum_{i=1}^p\left\{
\begin{array}{l}
2\langle R_{e_i,Y}B(X,e_i),W\rangle +2\langle R_{e_i,X}B(Y,e_i),W\rangle\\
-\langle A^W(X),R_{e_i,Y}e_i\rangle- \langle A^W(Y), R_{e_i,X}e_i\rangle\\
+\langle R_{e_i,B(X,Y)}e_i, W\rangle -2\langle A^W(e_i), R_{e_i,X}Y\rangle
\end{array}\right\}
\end{eqnarray*}
where $R$ is the curvature tensor for the ambient space. The principal result of this paper comes from a control of this term in a particular case. 

The case that we consider is where $X$ is a Riemannian symmetric space. We summarise some standard facts that we will use later. A fuller reference for this material is \cite{KN2}. A (connected) Riemannian manifold $X$ is a {\em Riemannian symmetric space} if for each point $p\in X$ there is an isometry $\sigma_p$ of $X$ that fixes $p$ and has derivative $-Id$ at $p$. This in particular implies that $X=G/H$ is homogeneous and the symmetry $\sigma_p$ induces an involutive automorphism $\sigma$ of the isometry group $G$ and hence of $\mathfrak{g}$. The Lie algebra $\mathfrak{g}$ splits 
\begin{eqnarray*}
\mathfrak{g}=\mathfrak{h}+\mathfrak{m}
\end{eqnarray*}
into the $+1$ and $-1$ eigenspaces of the automorphism. The ensemble $(\mathfrak{g},\mathfrak{h},\sigma)$ is referred to as a symmetric Lie algebra. The space $\mathfrak{m}$ can be identified with the tangent space to $X$ at a fixed point $p$. 

The algebraic structures of $(\mathfrak{g},\mathfrak{h},\sigma)$ can be related to the Riemannian geometry of $X$ by making the fundamental observation (see \cite{KN2}) that the set of tensors on $\mathfrak{m}$ that are invariant under the action of $\mathfrak{h}$ are in a one-to-one correspondence with the set of tensor fields on $X$ that are invariant under the action of $G$. 

We assume that $\mathfrak{g}$ is semi-simple, and that the Killing form is negative definite. This is the condition for the symmetric space to be of \emph{compact type}. The negative of the Killing form, restricted to $\mathfrak{m}$, defines a positive definite $\mathfrak{h}$-invariant bilinear form. It therefore corresponds to a $G$-invariant Riemannian metric on the space $X$. We will take this as our background metric. The Riemannian curvature of $X$ is given by, using the identification of $T_pX$ and $\mathfrak{m}$,
\begin{eqnarray*}
R_{X,Y}Z=-[[X,Y,]Z]
\end{eqnarray*}
The Ricci curvature is given by 
\begin{eqnarray*}
\text{Ric}(X,Y)=\tr_\mathfrak{m}(Z\mapsto -[X,[Y,Z]]).
\end{eqnarray*}
This is an $\mathfrak{h}$-invariant bilinear form on $\mathfrak{m}$. If $\mathfrak{h}$ acts irreducibly on $\mathfrak{m}$ this must be a multiple of the metric (as can be seen by simultaneously diagonalizing this with the metric). That is, $\text{Ric}(X,Y)=\rho\langle X,Y\rangle$ and necessarily $\rho>0$. In general, if $\mathfrak{g}$ is semi-simple, $\mathfrak{m}$ splits into the sum of irreducible representations of $\mathfrak{h}$ and the Ricci tensor is a multiple of the metric when restricted to each irreducible factor. That is,
\begin{eqnarray*}
\mathfrak{m}&=&\mathfrak{m}_1+\cdots +\mathfrak{m}_k,\\
\text{Ric}&=&\rho_1\langle\, ,\,\rangle|_{\mathfrak{m}_1}+\cdots+\rho_k\langle\, ,\,\rangle|_{\mathfrak{m}_k}
\end{eqnarray*}
and $\text{Ric}(X,Y)\geq \rho\langle X,Y\rangle$ for all $X,Y\in\mathfrak{m}$ where $\rho=\min_i\rho_i$. $\rho$ can be taken to be the smallest Ricci curvature of any unit tangent vector to the manifold $X$.

The {\em rank} of a symmetric space is defined to be the dimension of a maximal subspace $V\subseteq\mathfrak{m}$ for which $[X,Y]=0$ for all $X,Y\in V$. This is analogous to the dimension of a maximal torus in a Lie group. This also corresponds to the maximal dimension of a tangent subspace on which the sectional curvature vanishes identically. Accordingly, compact-type symmetric spaces of rank-one have strictly positive sectional curvature.

Furthermore, the compact-type rank-one spaces can be easily classified (see \cite{Chavel}). There are only $S^n$, $\RP^n$, $\CP^n$, $\HP^n$ and $\OP^2$. The final space is the $16$-dimensional Cayley plane. In contrast to the other examples it does not exist in an infinite family and has as set of isometries the exceptional Lie group $F_4$. We will from this point only consider the non-exceptional rank-one spaces. The important geometric property that these spaces have is that they admit the Hopf fibrations of spheres fibred by totally geodesic spheres. For example, one can define the map $S^{4n+3}\to \HP^n$ by sending a point to the quaternion line that it spans. The fibres of the map are of the form $S^{4n+3}\cap \Ham e$ for $e\in S^{4n+3}\subseteq \Ham^{n+1}$ and so are totally geodesic. 





\section{Rigidity of Sphere Factors}

We now consider the ambient space to be $X=S^p\times X_2$ where $X_2$ is a symmetric space of compact type and $M\subseteq X$ is a $p$-dimensional minimal submanifold. As above, we denote by $\rho$ the smallest Ricci curvature of any direction tangent to $X$. 

We also consider the projection $\pi_2:M\to X_2$ of the submanifold to the second factor. We also denote the derivative of this map by $\pi_2$. This acts as a bundle map on $M$, from $T_M$ to $T_{X_2}|_M$.  We will assume that the uniform size of this map is small in operator norm.

\begin{thm}\label{thm-tgisol}
There exists $\Lambda>0$ such that if $M$ is a $p$-dimensional closed minimal submanifold of $S^p\times X_2$ that satisfies 
\begin{eqnarray*}
\|\pi_2\|<\Lambda\\
\|A\|^2<\frac{\rho}{q}
\end{eqnarray*}
then $M=S^p\times \{x\}$ for some $x\in X_2$. Here again $q=2-\frac{1}{\text{dim} X_2}$.
\end{thm}
 

For the proof of this theorem we consider the second order equation of Simons given in the previous section. In the case at hand the term $R^\prime$ vanishes because the ambient space is symmetric. We can also control the term $R(A)$.

\begin{thm}\label{RA1forSS}
There exists $C=C(p, X_2)>0$ such that for any $p$-dimensional closed minimal submanifold $M$ of the symmetric space $\overline{M}=S^p\times X_2$  for which $\|\pi_2\|\leq\Lambda$ the term $\overline{R}(A)$ satisfies
\begin{eqnarray*}
\langle R(A), A\rangle \geq \Big( (2\rho+\frac{1}{p-1})-C\Lambda^2\Big)\|A\|^2.
\end{eqnarray*}
\end{thm}

This is the main technical result of this paper and the proof will be given in Section \ref{Sec:proof}.

\proof{(Of Theorem \ref{thm-tgisol}) The proof is quite short. Take $\Lambda^2 =(\rho +\frac{1}{p-1})/C$.
 In this case, if $\|\pi_2\|\leq \Lambda$ we have $\langle \overline{R}(A), A\rangle \geq \rho \|A\|^2$.
 We consider the equation of Simons and using integration by parts we see
\begin{eqnarray*}
0\leq  \int_M\|\nabla A\|^2 &\leq& -\int_M \langle A, \nabla^2A\rangle\\
&=& \int \langle A\circ \tilde{A} + \utA\circ A -\overline{R}(A), A\rangle\\
&\leq & \int q\|A\|^4 -\rho\|A\|^2\\
&=& q\int \|A\|^2(\|A\|^2 -\frac{\rho}{q}).
\end{eqnarray*}
Thus, if $\|A\|^2<\rho/q$ uniformly on $M$ we must have $A\equiv 0$ and $M$ is totally geodesic. That is, for this value of $\Lambda$ the hypotheses imply that $M$ is totally geodesic. By theorem \ref{T:tgrigid}
 we can take a $\Lambda$ so that we can conclude that $M=S^p\times\{pt\}$
 }
 
A similar theorem can be given  where we consider the intrinsic scalar curvature $K$ of $M$. We note that this curvature is with repsect to the metric induced on $M$ from the ambient space. 

\begin{thm}\label{T:Kest-sphere}
There exists $\Lambda>0$ such that if $M^p\subseteq S^p\times X_2$ is a closed minimal submanifold that satisfies
\begin{eqnarray*}
\|\pi_2\|&<&\Lambda\\
\frac{p}{2}-K&<&\frac{\rho}{q}
\end{eqnarray*}
then $M=S^p\times\{pt\}$.
\end{thm}

\proof{If one inspects the proof of Theorem \ref{thm-tgisol}, one can see that, for a given $\Lambda>0$, if $\|\pi_2\|<\Lambda$ and 
\begin{eqnarray*}
\|A\|^2-\frac{(2\rho +\frac{1}{p-1})-C\Lambda^2}{q}<0
\end{eqnarray*}
uniformly on $M$, then necessarily $A\equiv 0$. By the Gauss equation, we have
\begin{eqnarray*}
\|A\|^2=\sum_{ij}\|B(e_i,e_j)\|^2&=&\sum_{i\neq j}\langle R^X_{e_i,e_j}e_j,e_i\rangle-\sum_{i\neq j}\langle R^M_{e_i,e_j}e_j,e_i\rangle\\
&=&\sum_{i\neq j}\langle R^X_{e_i,e_j}e_j,e_i\rangle - K.
\end{eqnarray*}
 The ambient curvatures are obtained from the two factors by
\begin{eqnarray*}
\langle R^X_{e_i,e_j}e_j,e_i\rangle &=& \|[\pi_1e_i,\pi_1e_j]\|^2+\|[\pi_2e_i,\pi_2e_j]\|^2\\
&\leq & \frac{1}{2(p-1)} +\Lambda^2K_2^2.
\end{eqnarray*} 
Here $K_2$ is the maximum sectional curvature of a plane tangent to $X_2$ (see Section \ref{Sec:proof}). Thus, 
\begin{eqnarray*}
\|A\|^2 &\leq & \frac{p}{2} +p(p-1)K^2_2\Lambda^2 -K,\\
\|A\|^2 -\frac{(2\rho +\frac{1}{p-1})-C\Lambda^2}{q}&\leq& \frac{p}{2} -K-\frac{(2\rho +\frac{1}{p-1})-(C+p(p-1)K_2^2)\Lambda^2}{q}.
\end{eqnarray*}
We can take $\Lambda>0$ so that the right hand side equals 
\begin{eqnarray*}
\frac{p}{2} -K-\frac{\rho}{q}.
\end{eqnarray*}
Thus, if $\frac{p}{2}-K<\rho/q$ and $\|\pi_2\|<\Lambda$ uniformly on the submanifold, then necessarily $A\equiv 0$. As before, by Theorem \ref{T:tgrigid}, we can find a $\Lambda>0$ so that $M$ must be $S^p\times \{pt\}$.
}


\section{Isolation of Totally Geodesic Factors in Products}

In this section we consider the totally geodesic submanifolds of symmetric spaces. It is a basic result of the theory that complete totally geodesic submanifolds of a symmetric space $X$ with corresponding symmetric Lie algebra $\mathfrak{g}=\mathfrak{m}+\mathfrak{h}$ are in a one-to-one correspondence with subspaces $\mathfrak{t}\subseteq \mathfrak{m}$ that satisfy $[[\mathfrak{t},\mathfrak{t}],\mathfrak{t}]\subseteq\mathfrak{t}$. These subspaces are called \emph{Lie Triple Systems}.

Let $(\mathfrak{g},\mathfrak{h},\sigma)$ be a symmetric Lie algebra with $\mathfrak{g}$ semi-simple and of compact type. We will assume that this structure is reducible in that $\mathfrak{g}$ splits as 
\begin{eqnarray*}
\mathfrak{g}=(\mathfrak{m}_1+\mathfrak{h}_1)+(\mathfrak{m}_2+\mathfrak{h}_2).
\end{eqnarray*}
We consider the case that the first factor is the symmetric Lie algebra for the round sphere. That is, 
$(\mathfrak{g}_1,\mathfrak{h}_1,\sigma)=(\mathfrak{so}(n+1),\mathfrak{so}(n),\sigma)$. The important features of this space are that it has rank one and that we know all of its Lie triple systems. The totally geodesic subspaces of the sphere are the great spheres of the various dimensions and so have corresponding symmetric Lie algebras (conjugate to) $(\mathfrak{so}(p+1),\mathfrak{so}(p),\sigma)$. For notational reasons we will continue to refer to the first factor as $\mathfrak{m}_1$.

\begin{prop} Let $\mathfrak{t}$ be a Lie triple system contained in $\mathfrak{m}_1+\mathfrak{m}_2$. Consider the orthogonal projection $\pi_2:\mathfrak{t}\to \mathfrak{m}_2$. Suppose that $\|\pi_2\|\leq \Lambda<1$. 

Then the subalgebra $\mathfrak{t}+[\mathfrak{t},\mathfrak{t}]$ is simple and isomorphic to $\mathfrak{so}(p+1)$.
\end{prop}

\proof{
It is clear that the subspace $\mathfrak{k}= \mathfrak{t}+[\mathfrak{t},\mathfrak{t}] $ is a subalgebra of $\mathfrak{g}$. We consider the totally geodesic subspace $\pi_1(\mathfrak{t})\subseteq \mathfrak{m}_1$. This is clearly a Lie triple system of $\mathfrak{m}_1$ and so corresponds to a totally geodesic subspace of $X_1$.
 By hypothesis, $X_1=S^n$ and the corresponding totally geodesic subspace is a great sphere. This has isometry group $SO(p+1)$. 
 
 We consider $\pi_1:\mathfrak{t}\to\mathfrak{m}_1$. The assumption that $\|\pi_2\|< 1$ implies that $\pi_1$ is injective. We claim that $\pi_1$ is injective when considered on $\mathfrak{k}=\mathfrak{t}+[\mathfrak{t},\mathfrak{t}]$.  Let $x=X_1+X_2,\, y=Y_1+Y_2\in \mathfrak{t}$ where $X_i,\,Y_i\in\mathfrak{m}_i$. We can suppose that $X_1$ and $Y_1$ are non-zero. Then,
 \begin{eqnarray*}
[x,y]=[X_1,Y_1]+[X_2,Y_2]\\
\text{ and }  \pi_1[x,y]=[X_1,Y_1].
\end{eqnarray*}
Suppose that $[X_1,Y_1]=0$. We have assumed that the symmetric Lie algebra $(\mathfrak{g},\mathfrak{h}_1,\sigma)$ has rank one. This means that the dimension of a maximal subspace of $\mathfrak{m}_1$ on which the brackets vanish is equal to one. In other words,
\begin{eqnarray*} 
[X_1,Y_1]=0 \implies X_1=\lambda Y_1.
\end{eqnarray*}
We can rescale $x$ and $y$ so that $\lambda=1$. Then, $x-y=X_2-Y_2\in \text{ker}\pi_1\cap\mathfrak{t}$. This implies that $x=y$ and $[x,y]=0$. That is, $\pi_1$ is injective on $[\mathfrak{t},\mathfrak{t}]$ and so $\pi_1:\mathfrak{k}\to\mathfrak{m}_1+\mathfrak{h}_1$ is an isomorphism to its image. 

Hence, $\mathfrak{k}\cong \mathfrak{so}(p+1)$ and the symmetries correspond.
}

In particular the algebra $\mathfrak{k}$ is simple so the map $\pi_2:\mathfrak{k}\to\mathfrak{m}_2+\mathfrak{h}_2$ is either identically zero or an isomorphism to its image. In the first case, if $T\subseteq X_1\times X_2$ is the corresponding totally geodesic subspace, $\pi_2(T)$ is a point. In the second case, $\pi_2:T\to \pi_2(T)$ is a covering map. We show that if $\|\pi_2\|<\Lambda$ for $\Lambda$ sufficiently small the second case cannot occur.

\begin{thm}\label{T:tgrigid}
There exists $\Lambda>0$ such that if $T$ is a $p$-dimensional totally geodesic submanifold of $S^n\times M_2$ and $\|\pi_2\|\leq\Lambda$ then $T\subseteq S^n\times\{pt\}$. 
\end{thm}

\proof{
We suppose that $\pi_2(T)$ is not a point. In this case, it is a totally geodesic submanifold of the same dimension as $T$. By the area formula \cite{Simon} we have
\begin{eqnarray*}
\text{vol}(\pi_2(T))=\int_{\pi_2(T)}d\mathcal{H}^p(y)\leq \int_{\pi_2(T)}\int_{\pi_2^{-1}(y)}d\mathcal{H}^0(t)d\mathcal{H}^p(y)=\int_{T}J_{\pi_2}(x)d\mathcal{H}^p(x).
\end{eqnarray*}
The Jacobian terms are given, in the current case by,
\begin{eqnarray*}
(J_{\pi_2})^2&=&\det(d\pi_2))^*(d\pi_2))\leq \Lambda^{2p}\\
\text{ and }\ \ \text{vol}(\pi_2(T))&\leq& \Lambda^p\text{vol}(T).
\end{eqnarray*}
Similarly, by considering the projection to the other factor one can see that 
\begin{eqnarray*}
\text{vol}(T)\leq \frac{1}{(1-\Lambda^2)^{1/2}}\text{vol}(S^p)
\end{eqnarray*}
where $S^p$ has the metric induced from that on $S^n$ in this case. Thus,
\begin{eqnarray}\label{E:tgest}
\text{vol}(\pi_2(T))\leq \Big( \frac{\Lambda^2}{1-\Lambda^2}\Big)^\frac{p}{2}\text{vol}(S^p).
\end{eqnarray}
However, we can note that for (complete) totally geodesic submanifolds, the ambient geodesic spray from a tangent plane maps to the submanifold. As such, we can also note that the function $F$ defined by 
\begin{eqnarray*}
F(V)=\text{vol} \Big(\exp_o(U\cap V)\Big), 
\end{eqnarray*}
is continuous and bounded away from zero, where $V\in G(p,T_oX_2)$ and $U$ is a fixed open set containing the origin in $T_oX_2$ containing no tangential cut points. This fact, together with Equation \ref{E:tgest} implies that if $T\subseteq S^n\times X_2$ is totally geodesic and $\|\pi_2\|<\Lambda$ for sufficiently small $\Lambda$, then $\pi_2\equiv 0$. 

}


\section{Proof of Theorem \ref{RA1forSS}}\label{Sec:proof}

In this section we give the proof of Theorem \ref{RA1forSS}. We restate it here.

\begin{defn} 
Let $X=X_1\times X_2$ be a compact-type Riemannian symmetric space with metric induced from the Killing form. Let $\mathfrak{g}=\mathfrak{h}_1+\mathfrak{m}_1+\mathfrak{h}_2+\mathfrak{m}_2$ denote the decomposition of the Lie algebra of Killing fields. Define 
\begin{eqnarray*}
K_1&=&\max \{\|[X,Y]\|;\,X,\,Y\in\mathfrak{m}_1,\ |X|=|Y|=1\}\\
&=&\max \{\text{sec}(\sigma);\, \sigma \text{ is a plane tangent to }X_1\}\\
K_2&=&\max \{\|[X,Y]\|;\,X,\,Y\in\mathfrak{m}_2,\ |X|=|Y|=1\}
\end{eqnarray*}
\end{defn}

\begin{thm3}
There exists $C=C(p, M_2)>0$ such that for any $0<\Lambda\leq 1$ and for any $p$-dimensional closed minimal submanifold $M$ of the symmetric space $\overline{M}=S^p\times M_2$  for which $\|\pi_2\|\leq\Lambda$ the term $R(A)$ satisfies
\begin{eqnarray*}
\langle R(A), A\rangle \geq \Big( (2\rho+\frac{1}{p-1})-C\Lambda^2\Big)\|A\|^2.
\end{eqnarray*}
\end{thm3}

The term $R(A)$ is an section of the bundle $\Hom(N_M,S(M))$ and is given by the expression
\begin{eqnarray}\label{E:RAdefn}
\langle R(A)^WX,Y\rangle &=& 
\sum_{i=1}^p\left\{
\begin{array}{l}
2\langle R_{e_i,Y}B(X,e_i),W\rangle +2\langle R_{e_i,X}B(Y,e_i),W\rangle\\
-\langle A^W(X),R_{e_i,Y}e_i\rangle- \langle A^W(Y),R_{e_i,X}e_i\rangle\\
+\langle R_{e_i,B(X,Y)}e_i, W\rangle -2\langle A^W(e_i), R_{e_i,X}Y\rangle
\end{array}\right\}
\end{eqnarray}
That is, $R(A)=(1)+\cdots+(6)$. We will calculate the inner products $\langle(1),A\rangle$ and compare them in each case to $\|A\|^2$. For example,
\begin{eqnarray*}
\langle (1)^WX,Y\rangle &=& 2\sum_i\langle R_{e_i,Y}B(X,e_i),W\rangle
\end{eqnarray*}

The first observation that we make on these terms is of the symmetry between some of them.
\begin{lemma}
\begin{eqnarray*}
\langle (2),A\rangle &=& \langle (1),A\rangle\\
\langle (4),A\rangle &=& \langle (3),A\rangle.
\end{eqnarray*}
\end{lemma}
\proof{ This follows immediately by observing that in Equation \ref{E:RAdefn} , for fixed normal vector $W$, the terms $(1)$ and $(2)$, and $(3)$ and $(4)$ are respectively transposes of one another. They will then have the same inner product with the symmetric transformation $A$. 
}

We let $N$ be the dimension of $X$.
\begin{lemma}\label{lemma-1est}
Let $M$ be a $p$-dimensional minimal submanifold of the symmetric space $S^p\times X_2$. Suppose that $\pi_2$ satisfies $\|\pi_2\|\leq \Lambda$. Then the second fundamental form for $M$ satisfies 
\begin{eqnarray*}
\langle (1), A\rangle &=& \geq -2p^2(N-p)(K_1^2+K_2^2)\Lambda^2\|A\|^2.
\end{eqnarray*}
\end{lemma}

\proof{The term $(1)$ is defined by
\begin{eqnarray*}
\langle (1)^WX,Y\rangle &=& -2\sum_i\langle [[e_i,Y], B(X,e_i)],W\rangle\\
&=& 2\sum_i\langle [e_i, [B(X,e_i),W]],Y\rangle\\
\langle (1), A\rangle &=& \sum_{j,k}\langle (1)^{\eta_j}e_k,A^{\eta_j}(e_k)\rangle\\
&=& 2\sum_{ijk}\langle [A^{\eta_j}(e_k), e_i],[B(e_i,e_k),\eta_j]\rangle.
\end{eqnarray*}
We now note that if the symmetric Lie algebra splits as $\mathfrak{g}=\mathfrak{h}_1+\mathfrak{m}_1+\mathfrak{h}_2+\mathfrak{m}_2$ the terms in the Lie bracket calculation above are given as the sum from the respective factors. That is,
\begin{eqnarray*}
\langle (1),A\rangle &=&2\sum_{ijk}\langle [\pi_1A^{\eta_j}(e_k), \pi_1e_i],[\pi_1B(e_i,e_k),\pi_1\eta_j]\rangle\\
&&\ \ \ \ +2\sum_{ijk}\langle [\pi_2A^{\eta_j}(e_k), \pi_2e_i],[\pi_2B(e_i,e_k),\pi_2\eta_j]\rangle\\
&=& \langle (1),A\rangle_1+\langle(1),A\rangle_2.
\end{eqnarray*}
We assume that the projection $\pi_2$ defined on the tangent space satisfies $\|\pi_2\|\leq \Lambda$. One can note that if $M$ is the same dimension as the first factor, this implies that $\pi_1$, when acting on the normal bundle, also has norm bounded by $\Lambda$. Then,
\begin{eqnarray*}
\langle (1),A\rangle_1&=&2\sum_{ijk}\langle [\pi_1A^{\eta_j}(e_k), \pi_1e_i],[\pi_1B(e_i,e_k),\pi_1\eta_j]\rangle\\
&\geq & -2\sum_{ijk} K_1|\pi_1A^{\eta_j}(e_k)||\pi_1e_i|\cdot K_1 |\pi_1B(e_i,e_k)||\pi_1\eta_j|\\
&\geq& -2p^2(N-p) K_1^2 \Lambda^2\|A\|^2
\end{eqnarray*}
An identical calculation is made for the second term.
}

\begin{lemma}
If $M$ is a $p$-dimensional minimal submanifold of $S^p\times X_2$ then the second fundamental form of $M$ satisfies
\begin{eqnarray*}
\langle (3),A\rangle &\geq& \rho \|A\|^2-p(N-p)^2(K_1^2+K_2^2)\Lambda^2\|A\|^2.
\end{eqnarray*}
\end{lemma}
\proof{$(3)$ is give by
\begin{eqnarray*}
\langle (3),A\rangle &=&-\sum_i \langle A^W(X),R_{e_i,Y}e_i\rangle.
\end{eqnarray*}
The term on the right hand side looks very much like the ambient Ricci curvature operator, except that one must note that the trace is only over the tangent space to the submanifold rather than the ambient space. We can get around this fact by assuming that the submanifold is close to an irreducible factor.

\begin{eqnarray*}
\langle (3)^WX,Y\rangle &=& \sum_i \langle A^W(X), R_{e_i,Y}e_i\rangle\\
&=& -\sum_i\langle [[A^W(X),e_i],e_i],Y\rangle\\
\langle (3),A\rangle &=& -\sum_{ijk} \langle [[A^{\eta_j}(e_k),e_i],e_i],A^{\eta_j}(e_k)\rangle\\
&=& -\sum_{ijk}\langle \ad(A^{\eta_j}(e_k))\circ\ad(A^{\eta_j}(e_k))e_i,e_i\rangle\\
&=&-\sum_{jk}\tr_\mathfrak{m}((\ad(A^{\eta_j}(e_k))^2)-\sum_{jkl}\|[A^{\eta_j}(e_k),\eta_l]\|^2\\
&&\text{since $\{e_i,\, \eta_l\}$ forms an orthonormal basis for $T_X$,}\\
&=&\sum_{jk} \text{Ric}(A^{\eta_j}(e_k), A^{\eta_j}(e_k))-\sum_{jkl}\|[A^{\eta_j}(e_k),\eta_l]\|^2\\
&\geq& \rho \|A\|^2 -\sum_{jkl}\|[A^{\eta_j}(e_k),\eta_l]\|^2
\end{eqnarray*}
As in the proof of Lemma \ref{lemma-1est}, one can estimate the remaining term.
\begin{eqnarray*}
\|[A^{\eta_j}(e_k),\eta_l]\|^2&=&\|[\pi_1 A^{\eta_j}(e_k),\pi_1\eta_l]\|^2+\|[\pi_2A^{\eta_j}(e_k),\pi_2\eta_l]\|^2\\
&\leq& \Lambda^2 (K_1^2+K_2^2)\|A^{\eta_j}(e_k)\|^2\leq (K_1^2+K_2^2)\Lambda^2\|A\|^2\\
\text{so }\ \ \ \ \ \langle (3),A\rangle &\geq& \rho\|A\|^2-p(N-p)^2(K_1^2+K_2^2)\Lambda^2\|A\|^2.
\end{eqnarray*}
}
\begin{lemma}
The fifth factor of $R(A)$ satisfies 
\begin{eqnarray*}
\langle (5),A\rangle \geq -p^3(K_1+K_2^2)\Lambda^2\|A\|^2.
\end{eqnarray*}
\end{lemma}

\proof{ The fifth term in the expression for $R(A)$ is given by 
\begin{eqnarray*}
\langle (5)^WX,Y\rangle &=& \sum_i \langle R_{e_i,B(X,Y)}e_i, W\rangle\\
&=& -\sum_i\langle [[e_i, B(X,Y)],e_i], W\rangle\\
&=& \sum_i \langle A_X(\pi^N[e_i, [e_i,W]]),Y\rangle
\end{eqnarray*}
\begin{eqnarray*}
\langle (5),A\rangle& =&\sum_{ijk}\langle A_{e_k}(\pi^N[e_i,[e_i,\eta_j]]), A_{e_k}(\eta_j)\rangle\\
& =&\sum_{ijkl}\langle A_{e_k}(\pi^N[e_i,[e_i,\eta_j]]),e_l\rangle\langle e_l, A_{e_k}(\eta_j)\rangle\\
&=& \sum_{ijkl} \langle [e_i,[e_i,B(e_k,e_l)]], \eta_j\rangle \langle \eta_j, B(e_k,e_l)\rangle\\
&=& -\sum_{ikl} \|[B(e_k, e_l),e_i]\|^2
\end{eqnarray*}
The desired inequality can be obtained by observing
\begin{eqnarray*}
\|[B(e_k,e_l),e_i]\|^2&=&\|[\pi_1B(e_k,e_l),\pi_1e_i]\|^2+\|[\pi_2B(e_k,e_l),\pi_2e_i]\|^2\\
&\leq & \Lambda^2 \|B(e_k,e_l)\|^2K_1^2 +\Lambda^2 \|B(e_k,e_l)\|^2K_2^2\leq (K_1^2+K_2^2)\|A\|^2\Lambda^2 ,
\end{eqnarray*}
\begin{eqnarray*}
\text{and so }\ \ \ \ \langle (5), A\rangle\geq -p^3(K_1^2+K_2^2)\Lambda^2 \|A\|^2.
\end{eqnarray*}
}
Note here that the hypothesis $\|\pi_2\|\leq \Lambda$ implies that $|\pi_2(e)|\leq \Lambda$ for unit tangent vectors and $|\pi_1\eta|\leq\Lambda$ for unit normal vectors. This holds only if the dimension of the submanifold is the same as that of the first factor. We have used this in each of the above calculations.

Also note that to this point we have not used the fact that the first factor is the round sphere. We require this in the following calculation.  In this case we also consider the maps $\pi^T\pi_2$ and $\pi^T\pi_1$ where $\pi^T$ is the projection from $T_X$ to $T_M$. The assumption that $\|\pi_2\|\leq\Lambda$ is equivalent to the requirement that $\sum_i\lambda_i^2\leq \Lambda^2$ where the $\lambda^2_i$'s are eigenvalues of $\pi^T\pi_2$.

\begin{lemma}  \label{lemma-6est}
Under the assumption that $0\leq \Lambda<1$ uniformly on $M$, the sixth and final term of $\overline{R}(A)$ satisfies
\begin{eqnarray*}
\langle (6), A\rangle \geq \frac{1}{p-1} \|A\|^2-\Big(\frac{p^2+2}{p-1}+2p^2(N-p)K_2^2\Big)\Lambda^2\|A\|^2.
\end{eqnarray*}
\end{lemma}

\proof{
The term $(6)$ satisfies
\begin{eqnarray*}
\langle (6)^WX,Y\rangle &=& -2 \sum_i\langle A^W(e_i), R_{e_i,X}Y\rangle\\
&=& 2\langle A^W, [[e_i,X],Y]\rangle\\
&=& 2\sum_i\langle -[[e_i,X],A^W(e_i)],Y\rangle\\
\text{so }\ \ \ \ \langle (6),A\rangle &=&\langle (6),A\rangle_1+\langle (6),A\rangle_2
\end{eqnarray*}
That is,
\begin{eqnarray*}
\langle (6),A\rangle_1 &=& 2\sum_{ijk}\langle \frac{1}{2(p-1)}\big(-\langle \pi_1 e_i, \pi_1 A^{\eta_j}(e_i)\rangle\pi_1 e_k+\langle \pi_1e_k, \pi_1A^{\eta_j}(e_i)\rangle \pi_1e_i\big), \pi_1A^{\eta_j}(e_k)\rangle
\end{eqnarray*}
\begin{eqnarray*}
(p-1)\langle (6),A\rangle_1&=&
\sum_{ijk}\Large -\langle \pi^T\pi_1 e_i, A^{\eta_j}(e_i)\rangle \langle \pi^T\pi_1e_k, A^{\eta_j}(e_k)\rangle\\
&& \ \ \ \ \ \ +\langle B(\pi^T\pi_1 e_k, e_i), \eta_j\rangle \langle \eta_j, B(\pi^T\pi_1 e_i, e_k)\rangle\\
&=& -\sum_{ijk} \lambda_k^2\lambda_i^2\langle e_i, A^{\eta_j}(e_i)\rangle\langle e_k, A^{\eta_j}(e_k)\rangle\\
&&\ \ \ \ \ \ +\sum_{ik}(1-\lambda_i^2)(1-\lambda_k^2)\|B(e_i,e_k)\|^2\\
&\geq& -\sum_{ik}\lambda_i^2\lambda_k^2\langle B(e_i,e_i),B(e_k,e_k)\rangle +\sum_{ik}\|B(e_i,e_k)\|^2- \sum_{ik}(\lambda_i^2+\lambda_k^2)\|B(e_i,e_k)\|^2\\
&\geq& \|A\|^2-p^2\Lambda^2\|A\|^2-2\Lambda^2\|A\|^2.
\end{eqnarray*}
In the third equality we have used that the submanifold is minimal.

\begin{eqnarray*}
\langle (6),A\rangle_2 &=& -2\sum_{ijk} \langle [[\pi_2e_i, e_k], \pi_2 A^{\eta_j}(e_i)], A^{\eta_j}(e_k)\rangle\\
&=& -2 \sum_{ijk}\langle [\pi_2 e_i, e_k],[\pi_2 A^{\eta_j}(e_i), A^{\eta_j}(e_k)]\rangle\\
&\geq& -2\sum_{ijk} K_2\Lambda \cdot K_2\Lambda \|A^{\eta_j}(e_i)\| \|A^{\eta_j}(e_k)\|\\
&\geq & -2 p^2(N-p) K_2^2 \Lambda^2 \|A\|^2.
\end{eqnarray*}
}

Collecting the results from the Lemmas \ref{lemma-1est} to \ref{lemma-6est} if we define
\begin{eqnarray*}
C(p,X_2) =\Big( 4p^2(N-p)+2p(N-p)^2 +p^3\Big) \Large(K_1^2+K_2^2) +\Big(\frac{p^2+2}{p-1}+2p^2(N-p)K_2^2\Big),
\end{eqnarray*}
then we can conclude that the tensor $R(A)$ satisfies
\begin{eqnarray*}
\langle R(A),A\rangle \geq \Big( 2\rho +\frac{1}{p-1}- C\Lambda^2\Big)\|A\|^2.
\end{eqnarray*}
This concludes the proof of Theorem \ref{RA1forSS}.


\section{Riemannian Submersions}

In this section we review the basic definition of Riemannian submersions, with a particular interest in the way the tensors and curvatures behave for submanifolds of the spaces. We consider submersions $\pi:\overline{M}\to M$. That is, $\pi$ is surjective and has surjective derivative at all points in $\overline{M}$. The preimage of a point in $M$ is a submanifold of $\overline{M}$.

\begin{defn}
Let $\overline{M}$ and $M$ be Riemannian manifolds. The submersion $\pi$ is a \emph{Riemannian submersion} if $\pi_*$ is an isometry when restricted to the orthogonal space to the fibres.  That is,
\begin{eqnarray*}
\pi_*:T_pF^\perp\to T_{\pi(p)}M \ \ \ \text{is an isometry}
\end{eqnarray*}
for all fibres $F=\pi^{-1}(\pi(p))$ of the submersion.
\end{defn}
We call $T_pF$ the \emph{vertical} space and $T_pF^\perp$ the \emph{horizontal} space.  We respectively denote by $\mathcal{V}$ and $\mathcal{H}$ the projections to these spaces.

We assume that the fibres of the Riemannian submersion are totally geodesic.

\begin{defn} The \emph{O'Neill tensor} of the submersion is given by 
\begin{eqnarray*}
A_XY =\mathcal{V}\nabla_{\mathcal{H}X}\mathcal{H}Y + \mathcal{H}\nabla_{\mathcal{H}X}\mathcal{V}Y.
\end{eqnarray*}
\end{defn}

Let $\overline{R}$ and $R$ be the Riemannian curvatures of $\overline{M}$ and $M$ respectively. Also denote by 
$\overline{S}(X,Y)=\langle \overline{R}_{Y,X}X,Y\rangle$ for $X$ and $Y$ tangent to $\overline{M}$, and $S(X,Y)=\langle R_{Y,X}X,Y\rangle$ for $X$ and $Y$ tangent to $M$. Then, from \cite{Lawson2}, we have
\begin{eqnarray*}
S(\pi_*X, \pi_* Y) &=&\overline{S}(X,Y) +3|A_XY|^2,\\
S(X,V)&=& |A_XV|^2
\end{eqnarray*}
where $X$ and $Y$ are horizontal vectors and $V$ is vertical.

For $p\in\overline{M}$ we denote
\begin{eqnarray*} 
\overline{K}(p) &=& \text{ scalar curvature of $\overline{M}$ at $p$},\\
K(p)&=& \text{scalar curvature of $M$ at $\pi(p)$},\\
\tau (p)&=& \sum_{j,k}\overline{S}(e_j,\nu_k),\\
&=& \text{twisting curvature at $p$}\\
r(p)&=& \text{scalar curvature of the fibre $\pi^{-1}(\pi(p))$ at $p$.}
\end{eqnarray*}

\begin{thm}\emph{\cite{Lawson2}}\ \ \
$\ \ \ \ K=\overline{K}+\tau-r.$
\end{thm}
 We combine the study of submersions with that of submanifolds. In particular we suppose that $\overline{M}$ is a submanifold of $\overline{X}$ and $M$ is a submanifold of $X$ such that the following diagram commutes.  
\begin{eqnarray*}
\begin{CD}
\overline{M} @>{\bar f}>> \overline{X}\\
@V\pi VV   @VV\pi V\\
M @>f>> X
\end{CD}
\end{eqnarray*}
We suppose the fibres are totally geodesic in each case and that ${\bar f}$ is a diffeomorphism on the fibres.

\begin{thm} \emph{\cite{Lawson2}}
$M$ is a minimal submanifold of $X$ if and only if $\overline{M}$ is a minimal submanifold of $\overline{X}$. 
\end{thm}
We can use the Gauss equation to compare the twisting curvature for the submersions $\overline{M}\to M$ and $\overline{X}\to X$. We have
\begin{eqnarray*}
S_X(e_j, \nu_k)-S_M(e_j,\nu_k)=\|B(e_j,\nu_k)\|^2 - \langle B(e_j, e_j), B(\nu_k, \nu_k)\rangle=\|B(e_j, \nu_k)\|^2\geq 0.
\end{eqnarray*}
Finally, we note that if $\overline{M}\subseteq \overline{X}=\overline{X}_1\times X_2$ and $M\subseteq X=X_1\times X_2$ then the projections $\pi_2$ from $M$ and $\overline{M}$ respectively to $X_2$ have the same uniform norm.


\section{Rigidity of Rank-One Factors}

In this section we extend the previous calculations to the case where the submanifold $M$ is close to a rank-one factor of $X$ other than the sphere. That is, we take $X=X_1\times X_2$ where $X_1$ is one of either $\RP^n$, $\CP^n$ or $\HP^n$. The first case follows immediately, because the curvature tensor of $\RP^n$ coincides, up to a factor, with that of $S^n$. For the other spaces, we use the fact that they admit Hopf-fibrations of spheres with totally geodesic fibres. We can then use the Riemannian submersion framework developed in the previous section. 

\begin{lemma} \label{L:CPtauest}
Let $M$ be a $2n$-dimensional minimal submanifold of the compact-type symmetric space $\CP^n\times X_2$ and let $\overline{M}\subseteq S^{2n+1}\times X_2$ be the minimal submanifold that fibres over $M$. Then, 
\begin{eqnarray*}
\tau_M\leq \frac{1}{2}.
\end{eqnarray*}
\end{lemma}
\proof{
We can assume that $\{e_i\}$ and $\{\nu\}$ form bases for the horizontal and normal spaces respectively and   together diagonalize the map $\pi^T\pi_1$ on $\overline{M}$. In particular, $\pi_1\nu=\nu$ because $\nu$ is vertical and $|\pi_1 e_i|=\lambda_i$.
\begin{eqnarray*}
S_M(e_i,\nu)&\leq &S_X(e_i,\nu)\\
&=& \|[\pi_1e_i, \pi_1\nu]\|^2 +\|[\pi_2 e_i,\pi_2\nu]\|^2\\
&=& \frac{1}{4n}\lambda_i^2\leq \frac{1}{4n}.
\end{eqnarray*}
Thus, $\tau_M(p)\leq \frac{1}{4n}2n=\frac{1}{2}$.
}

\begin{thm}\label{CPnRigid}
There exists $\Lambda>0$ such that if $M$ is a $2n$-dimensional closed minimal submanifold of $\CP_n\times X_2$ that uniformly satisfies 
\begin{eqnarray*}
\|\pi_2\|&\leq &\Lambda\\
(n+1)-K^\prime &<&\frac{\rho}{q}.
\end{eqnarray*}
Then $M=\CP^n\times\{pt\}$ and $K \equiv n+1 $.
\end{thm}

\proof{Let $\overline{M}\subseteq S^{2n+1}\times X_2$ be the closed minimal submanifold that fibres over $M$. The fibres in this case are totally geodesic circles so $r=0$. Then,
\begin{eqnarray*}
\frac{2n+1}{2}-\overline{K} =\frac{2n+1}{2}-K+\tau\leq n+1-K
\end{eqnarray*}
so by Theorem \ref{T:Kest-sphere} we can conclude the result.
}

\begin{cor}
There is a $C^3$-open neighbourhood of the standard embedding of $\CP_n$ in $\CP_n\times M_2$ in the set of immersions such that any minimal immersion contained in it is conjugate to the standard one.
\end{cor}

We now consider closed $4n$-dimensional minimal submanifolds of $\HP^n\times X_2$. For such a minimal submanifold $M$ there is minimal submanifold $\overline{M}$ of $S^{4n+3}\times X_2$ that fibres over it. The fibres are totally geodesic in both $\overline{M}$ and $S^{4n+3}\times X_2$ and so, with respect to the metric that we have considered, have scalar curvature $\frac{3}{4n+2}$. In an identical way to Lemma \ref{L:CPtauest} we can calculate $\tau_M$.
\begin{lemma}For $p\in \overline{M}\subseteq S^{4n+3}\times X_2$,
\begin{eqnarray*}
r(p)&=&\frac{3}{4n+2}\\
\tau_M(p)&\leq& \frac{3n}{2n+1}.
\end{eqnarray*}
\end{lemma}

\begin{thm}\label{HPnrigid}
Let $M$ be a closed $4n$-dimensional minimal submanifold of $\HP_n\times M_2$ and let $K^\prime$ be its intrinsic scalar curvature. Suppose that $M$ satisfies, for $\Lambda$ given above, 
\begin{eqnarray*}
\|\pi_2\|\leq \Lambda,\\
\frac{4n(n+2)}{2n+1} -K^\prime <\frac{\rho}{q}.
\end{eqnarray*}
Then $M$ is totally geodesic. There exists a possibly smaller $\Lambda$ such that these hypotheses imply that $M$ is a totally geodesic factor $\HP_n\times pt$.
\end{thm}

The proof is identical to that of Theorem \ref{CPnRigid}.

\begin{cor}
There is a $C^3$-open neighbourhood of the standard embedding of $\CP_n$ in $\CP_n\times M_2$ in the set of immersions such that any minimal immersion contained in it is conjugate to the standard one.
\end{cor}

\bibliographystyle{amsalpha}

\begin{thebibliography}{EFK}


\bibitem[B] {Barbosa} J.L.M. Barbosa, {\em An extrinsic rigidity theorem for minimal immersions of $S^2$ into $S^n$}, J. Diff. Geom. {\bf 14} (1979) 355--368.


\bibitem[C] {Chavel} I. Chavel, {\em Riemannian Symmetric Spaces of Rank One}, Marcel Dekker Inc., New York 1972

\bibitem[CCK] {CCK} S.S. Chern, M. Do Carmo and S. Kobayashi, {\em Minimal submanifolds of the sphere with second fundamental form of constant length},  Functional Analysis and Related Fields, Springer Verlag, (1970) 59--75.


\bibitem[FC] {FC} D. Fischer-Colbrie, {\em Some rigidity theorems for minimal submanifolds of the sphere}, Acta Math. {\bf 145} (1980) 29--46. 

\bibitem[GMZ] {GMZ} H. Gluck, F. Morgan and W. Ziller, {\em Calibrated geometries in Grassmann manifolds}, Comment Math. Helv., {\bf 64} (1989) 256--268.

\bibitem[HH] {Hsiang} W.T. Hsiang and W.Y. Hsiang, {\em Examples of codimension-one closed minimal submanifolds in some symmetric spaces. I.} , J. Diff. Geom.  {\bf 15}  (1980) 543--551.

\bibitem[KN2] {KN2} S. Kobayashi and K. Nomizu, {\em Foundations of Differential Geometry}, volume 2. Interscience, New York, 1969.

\bibitem[L1]{Lawson1} H.B. Lawson Jr. {\em Complete minimal surfaces in $S^3$}, Ann. of Math., {\bf 92} (1970) 335--374. 

\bibitem[L2] {Lawson2} H.B. Lawson Jr., {\em Rigidity theorems in rank-1 symmetric spaces}, J. Diff. Geom., {\bf 4} (1970) 349-357.

\bibitem[S] {Simon} L. Simon, {\em Lecture Notes on Geometric Measure Theory}, Australian National University, 1983.

\bibitem[Si] {Simons} J. Simons, {\em Minimal varieties in riemannian manifolds}, Ann. Math., {\bf 88} (1968) 62--105.

\bibitem[T] {Thi} D.{\v C}. Thi, {\em Minimal real currents on compact riemannian manifolds}, Math. USSR Izv. {\bf 11} (1970) 807--820.



















\end{thebibliography}

\end{document}